\documentclass[12pt,amssymb,verbatim]{amsart}

\usepackage[all]{xy}

\def\cal{\mathcal}

\numberwithin{equation}{section}

\theoremstyle{plain}
\newtheorem{theorem}{Theorem}

\newtheorem{lemma}[theorem]{Lemma}

\newtheorem{corollary}[theorem]{Corollary}
\newtheorem{definition}[theorem]{Definition}

\theoremstyle{definition}
\newtheorem{question}[theorem]{Question}

\theoremstyle{remark}
\newtheorem{remark}[theorem]{Remark}

 \def\today{\ifcase\month\or
  January\or February\or March\or April\or May\or June\or
  July\or August\or September\or October\or November\or December\fi
  \space\number\day, \number\year}

\begin{document}
\title[Genus of division algebras]
{Genus of division algebras over \\ fields with infinite transcendence degree}
\author[Tikhonov]{Sergey V. Tikhonov}


\address{
Belarusian State University, Nezavisimosti Ave., 4,
220030, Minsk, Belarus} \email{tikhonovsv@bsu.by }

\def\cA{{\cal A}}
\def\cB{{\cal B}}
\def\cC{{\cal C}}
\def\cD{{\cal D}}
\def\cV{{\cal V}}
\def\cE{{\cal E}}
\def\cM{{\cal M}}
\def\cR{{\cal R}}
\def\cQ{{\cal Q}}
\def\M{{\rm M}}
\def\cS{{\cal S}}
\def\Symb{\textrm{Symb}}
\def\Gal{\textrm{Gal}}
\def\U{\textrm{U}}
\def\SU{\textrm{SU}}
\def\R{\textrm{R}}
\def\GL{\textrm{GL}}

\def\U{{\rm U}}
\def\SU{{\rm SU}}
\def\SL{{\rm SL}}
\def\op{{\rm op\;}}


\def\cor{\textrm{cor}}
\def\deg{\textrm{deg}}
\def\exp{\textrm{exp}}
\def\Gal{\textrm{Gal}}
\def\ram{\textrm{ram}}
\def\Spec{\textrm{Spec}}
\def\Proj{\textrm{Proj}}
\def\Perm{\textrm{Perm}}
\def\coker{\textrm{coker\,}}
\def\Hom{\textrm{Hom}}
\def\im{\textrm{im\,}}
\def\ind{\textrm{ind}}
\def\int{\textrm{int}}
\def\inv{\textrm{inv}}
\def\min{\textrm{min}}
\def\Br{\textrm{Br}}
\def\char{\textrm{char}}

\def\w{\widehat}

\begin{abstract}


We prove the finiteness of the genus of finite-dimensional division algebras over many infinitely generated fields.
More precisely, let $K$ be a finite field extension of a field which is a purely transcendental extension of infinite transcendence degree of some subfield.
We show that if $\cD$ is a central division $K$-algebra, then
${\bf gen}(\cD)$ consists of Brauer classes $[\cD']$ such that $[\cD]$ and $[\cD']$ generate the same subgroup of $\Br(K)$.
In particular, the genus of any division $K$-algebra of exponent 2 is trivial.
Note that the family of such fields is closed under finitely generated extensions.
Moreover, if $\char(K) \ne 2$, we prove that the genus of a simple algebraic group of type $\mathrm{G}_2$ over such a field $K$ is trivial.
\end{abstract}
\maketitle


\def\dd{{\partial}}

\def\into{{\hookrightarrow}}

\def\emptyset{{\varnothing}}

\def\alp{{\alpha}}  \def\bet{{\beta}} \def\gam{{\gamma}}
 \def\del{{\delta}}
\def\eps{{\varepsilon}}
\def\kap{{\kappa}}                   \def\Chi{\text{X}}
\def\lam{{\lambda}}
 \def\sig{{\sigma}}  \def\vphi{{\varphi}} \def\om{{\omega}}
\def\Gam{{\Gamma}}  \def\Del{{\Delta}}  \def\Sig{{\Sigma}}
\def\ups{{\upsilon}}


\def\A{{\mathbb A}}
\def\F{{\mathbb F}}
\def\Q{{\mathbb{Q}}}
\def\CC{{\mathbb{C}}}
\def\PP{{\mathbb P}}
\def\R{{\mathbb R}}
\def\Z{{\mathbb Z}}
\def\X{{\mathbb X}}

\def\Gm{{{\Bbb G}_m}}
\def\Gmk{{{\Bbb G}_{m,k}}}
\def\GmL{{\Bbb G_{{\rm m},L}}}
\def\Ga{{{\Bbb G}_a}}

\def\Fb{{\overline F}}
\def\Hb{{\overline H}}
\def\Kb{{\overline K}}
\def\Lb{{\overline L}}
\def\Yb{{\overline Y}}
\def\Xb{{\overline X}}
\def\Tb{{\overline T}}
\def\Bb{{\overline B}}
\def\Gb{{\overline G}}
\def\Vb{{\overline V}}

\def\kb{{\bar k}}
\def\xb{{\bar x}}

\def\Th{{\hat T}}
\def\Bh{{\hat B}}
\def\Gh{{\hat G}}

\def\Xt{{\tilde X}}
\def\Gt{{\tilde G}}

\def\gg{{\mathfrak g}}
\def\gm{{\mathfrak m}}
\def\gp{{\mathfrak p}}
\def\gq{{\mathfrak q}}

\def\textrm#1{\text{\rm #1}}

\def\res{\textrm{res}}
\def\cor{\textrm{cor}}
\def\R{\textrm{R}}
\def\ed{\textrm{ed}}

\def\tors{_{\textrm{tors}}}      \def\tor{^{\textrm{tor}}}
\def\red{^{\textrm{red}}}         \def\nt{^{\textrm{ssu}}}
\def\sc{^{\textrm{sc}}}
\def\sss{^{\textrm{ss}}}          \def\uu{^{\textrm{u}}}
\def\ad{^{\textrm{ad}}}           \def\mm{^{\textrm{m}}}
\def\tm{^\times}                  \def\mult{^{\textrm{mult}}}
\def\tt{^{\textrm{t}}}
\def\uss{^{\textrm{ssu}}}         \def\ssu{^{\textrm{ssu}}}
\def\cf{^{\textrm{cf}}}
\def\ab{_{\textrm{ab}}}

\def\et{_{\textrm{\'et}}}
\def\nr{_{\textrm{nr}}}

\def\SB{\textrm{SB}}

\def\<{\langle}
\def\>{\rangle}

\def\til{\;\widetilde{}\;}


\font\cyr=wncyr10 scaled \magstep1%
\def\Bcyr{\text{\cyr B}}
\def\Sh{\text{\cyr Sh}}
\def\Ch{\text{\cyr Ch}}

\def\lpsi{{{}_\psi}}
\def\bks{{\backslash}}

\def\Br{\textrm{Br}}
\def\Pic{\textrm{Pic}}
\def\Bt{{{}_2\textrm{Br}}}
\def\Bn{{{}_n\textrm{Br}}}
\def\Grass{\textrm{Grass}}






\section{Introduction} \label{sec:intro}

Let $F$ be a field and $\Br(F)$ its Brauer group. If classes of two central simple $F$-algebras generate the same subgroup of the group $\Br(F)$,
then these algebras have the same splitting fields.
Amitsur's theorem on generic splitting fields (\cite{Am55}) yields that if two central division $F$-algebras share the same splitting fields,
than these algebras have the same degree and their classes generate the same subgroup of $\Br(F)$.
However, if one considers only finite-dimensional splitting fields, the situation is not so clear.
In this regard, the notion of genus of a division algebra has been introduced.


\begin{definition}
The genus ${\bf gen}(\cD)$ of a finite-dimensional central division algebra $\cD$ over a field $F$ is defined as the set of classes $[\cD']\in \Br(F)$,
where $\cD'$ is a central division $F$-algebra having the same maximal subfields as $\cD$.
\end{definition}


This means that $\cD$ and $\cD'$ have the same degree $n$, and a field extension $K/F$ of degree $n$ admits an $F$-embedding $K \hookrightarrow \cD$ if and only if it admits an $F$-embedding
$K \hookrightarrow \cD'$. Different variations of the notion of the genus are mentioned in \cite{ChRaRa15} and \cite{KrMaRaRoSa22}.

If $F$ is a number field, then by Albert-Brauer-Hasse-Noether theorem (see \cite[\S 18.4]{Pi82}),
any central division $F$-algebra has finite genus and any quaterenion division $F$-algebra has trivial genus (see \cite[Proposition 3.1]{ChRaRa15}).
However, in the case of central division $F$-algebras of degree $\ge 3$ the genus might be arbitrary large (\cite[\S 1]{ChRaRa15}).

In \cite{ChRaRa13}, the authors describe a general approach to proving the finiteness of $\bold{gen}(D)$ and estimating its size that involves the unramified Brauer
group with respect to an appropriate set of discrete valuations of $K$.
In \cite{ChRaRa16} and \cite{ChRaRa20}, it is proved that the genus is finite for any central division algebra over a finitely generated field.

The genus of division algebras over discrete valued fields is considered in \cite{Sr24}. In particular, when the division
algebra is a quaternion algebra, the triviality of genus is proved over many fields which
include higher local fields, function fields of curves over higher local fields and
function fields of curves over real closed fields.

In \cite{RaRa10}, it is shown that if $|{\bf gen}(\cD)|=1$ for any central quaternion division algebra $\cD$ over a field $F$ of characteristic not 2,
then the same is true for central quaternion division algebras over the field $F(x)$.
The generalization of this result to central division algebras of exponent 2 is given in \cite{ChRaRa13}.

Thus the genus of a quaternion division algebra over a purely transcendental extension of a number field of any finite transcendence degree is trivial,
but the following question remains open.


\begin{question} (\cite[Question 8.8.]{RaRa20})
Does there exist a central quaternion division algebra $\cD$ over a finitely generated field of characteristic $\ne 2$ having nontrivial genus.
\end{question}


There is one more open question.


\begin{question} (\cite[Question 8.6.]{ChRaRa24}) \label{q:F(x)}
Does there exist a central quaternion division algebra $\cD$ over  the field of rational functions
$K=F(x)$ over some field $F$ having nontrivial genus ${\bf gen}(\cD)$.
\end{question}


Over general fields, the genus of algebras can be infinite.
In \cite{GaSa10}, the authors give an example of a quaternion algebra $\cD$ over a large center, constructed by iterative composition of function fields,
such that ${\bf gen}(\cD)$ does not consist of a single class.
In \cite{Me14}, it is shown that there exist quaternion division algebras with infinite genus over fields
with infinite transcendence degree over the prime subfield.
In \cite{Ti16}, these results are generalized to the case of division algebras of any prime degree.
In \cite[Theorem 5.2]{KrLi24}, the authors constructed examples of
infinite collections of division algebras over a field which share all the same finite
splitting fields and splitting fields which are function fields of curves of genus at
most 1.



In this paper, we prove the finiteness of the genus of finite-dimensional division algebras over many infinitely generated fields.
Note that the family of such fields is closed under finitely generated extensions.
The main result is the following


\begin{theorem} \label{th:main}
Let $K$ be a finite field extension of a field which is a purely transcendental extension of infinite transcendence degree of some subfield.
Let also $\cD$ be a central division $K$-algebra.
Then ${\bf gen}(\cD)$ consists of Brauer classes $[\cD']$ such that $[\cD]$ and $[\cD']$ generate the same subgroup of $\Br(K)$.

In particular, the genus of any division $K$-algebra of exponent 2 is trivial.
\end{theorem}



\begin{remark} \label{rem: genus}
If classes of central division $F$-algebras $\cD_1$ and $\cD_2$ generate the same subgroup of the Brauer group of the field $F$, then
${\bf gen}(\cD_1) = {\bf gen}(\cD_2)$. Thus in the case of central division algebras over fields from the previous theorem the genus is minimum possible.
\end{remark}



In a similar way one can define the genus of an absolutely almost simple algebraic group.

Let $F$ be a field and $F^{sep}$ its separable closure.
Two absolutely almost simple algebraic $F$-groups $G_1$ and $G_2$ are said to have the same $F$-isomorphism classes of
maximal $F$-tori if every maximal $F$-torus of $G_1$ is $F$-isomorphic to some maximal $F$-torus of $G_2$,
and vice versa. An algebraic $F$-group $G'$ is called a $F$-form of an algebraic $F$-group $G$ if $G$ and $G'$ are
isomorphic over  $F^{sep}$.

\begin{definition} \cite[Def. 6.1]{ChRaRa15}
Let $G$ be an absolutely almost simple algebraic group over a field $F$. The genus ${\bf gen}(G)$ of $G$ is the set of
$F$-isomorphism classes of $F$-forms $G'$ of $G$ that have the same $F$-isomorphism classes of maximal $F$-tori as $G$.
\end{definition}

The genus of an algebraic group is trivial in some special cases and it is conjectured to be finite whenever the field $F$ is finitely generated of "good" characteristic
(see details in \cite[\S 8]{RaRa20}).
Surveys of recent developments on the genus problem are given in \cite{RaRa20} and \cite{Ra22}.
Note that the genus can be infinite in the general case.
If $\cD$ is a finite-dimensional central division $F$-algebra, then it is well-known that any maximal $F$-torus of the corresponding algebraic group $G=\SL_{1,\cD}$
is of the form $\R_{E/F}(\mathbb{G}_m) \cap G$ (where $\R_{E/R}(\mathbb{G}_m)$ is the Weil restriction of the 1-dimensional split torus $\mathbb{G}_m$) for some maximal separable subfield $E$ of $\cD$.
Thus the results on genus of division algebras from \cite{Me14} and \cite{Ti16}
rephrased in the language of algebraic groups say that for any prime $p$, there exist fields
over which there are inner forms of type $\mathrm{A}_{p-1}$ with infinite genus.
In \cite{Ti24}, such an example is constructed for outer forms of type $\mathrm{A}_2$.
An example of groups of type $\mathrm{G}_2$ with infinite genus is obtained in \cite[Rem. 3.6(b)]{BeGiLe16}.
Note that fields in these examples have infinite transcendence degree over the prime subfield. In \cite[Theorem 1.6]{ChRaRa24}, it is shown that
if $G$ is a group of type $\mathrm{G}_2$ over a finitely generated field $F$ of characteristic $\ne 2$, then ${\bf gen}(G\times_F F(x_1,\dots, x_6 ))$  is trivial,
where $F(x_1,\dots, x_6)$ is the field of rational functions in 6 variables. It is also proved that the genus of a group of type $\mathrm{G}_2$ is trivial in case of some finitely generated fields
(see \cite[Theorem 8.14]{RaRa20}).

We obtained a new case of triviality of the genus of groups of type $\mathrm{G}_2$.
\begin{theorem} \label{th:G2}
Let $K$ be a finite field extension of a field which is a purely transcendental extension of infinite transcendence degree of some subfield $F$.
Assume $\char(K) \ne 2$. Let $G$ be a simple algebraic group of type $\mathrm{G}_2$ over $K$.
Then ${\bf gen}(G)$ is trivial.
\end{theorem}


Below for a field extension $E/F$ and a central simple $F$-algebra $\cA$, $\cA_E$ denotes the tensor product $\cA\otimes_F E$.




\section{Proof of main theorem} \label{sec:proof}

We start with the following

\begin{lemma} \label{l:finite subset}
Let $K$ be a finite field extension of a field which is a purely transcendental extension of infinite transcendence degree of some subfield.
Assume that $a_1,\dots,a_m \in K$.
Then there exists a subfield
$M \subset K$ such that
$a_1,\dots,a_m \in M$ and $K/M$ is a purely transcendental field extension with infinite transcendence degree.
\end{lemma}

\noindent {\it {Proof}}.
Consider a tower of field extensions
$$
F\subset E \subset K,
$$
were $E/F$ is a purely transcendental extension of infinite transcendence degree and
$K/E$ is a finite extension.

Let $S$ be a transcendence base of $E$ over $F$.
Let also $K=E(b_1,\dots,b_l)$, where $b_1,\dots,b_l$ are algebraic over $E$.
Then there exists a finite number of elements $s_1,\dots,s_k \in S$ such that
the elements $b_1,\dots,b_l,a_1,\dots,a_m$ are algebraic over the field  $F(s_1,\dots,s_k)$.
Let
$$
M:=F(s_1,\dots,s_k)(b_1,\dots,b_l,a_1,\dots,a_m).
$$
Since algebraically independent elements remain algebraically independent over algebraic extensions,
then $K/M$ is a purely transcendental field extension with
a transcendence base
$$
S \setminus \{ s_1,\dots,s_k \}.
$$
\qed


Now we are in a position to prove the main theorem.

%

\noindent {\it Proof of Theorem \ref{th:main}}.
Let $n$ be the degree of the algebra $\cD$
and $\cD'$ a central division $K$-algebra of degree $n$. 
Let also  $C$ be the finite set consisting of structure constants of the algebras $\cD$ and $\cD'$ relative to some bases of these algebras.
Be Lemma \ref{l:finite subset}, there exist a subfield
$M \subset K$ such that
$C \subset M$ and $K/M$ is a purely transcendental field extension with infinite transcendence degree.
Note that $\cD\cong \cB_K$ and $\cD'\cong \cB'_K$ for some central division $M$-algebras $\cB$ and $\cB'$.

Let $S$ be a transcendence base of $K$ over $M$. Let also $s_1,\dots,s_{n-1}$ be some elements from $S$
and
$$
L:= M (S \setminus \{ s_1,\dots,s_{n-1} \}).
$$
Then $K=L(s_1,\dots,s_{n-1})$.

Assume that the algebras $\cD$ and $\cD'$ have the same genus. Since $\cD\cong \cB_L \otimes_L K$ and $\cD'\cong \cB'_L \otimes_L K$,
as in the proof of Theorem 8.1 in \cite{ChRaRa24}, we obtain that
$[\cB_L]$ and $[\cB'_L]$ generate the same subgroup of $\Br(L)$.
Note that this part of the proof of Theorem 8.1 in \cite{ChRaRa24} does not use that the degree of the algebra is prime to the characteristic of the field.
Then  $[\cD]$ and $[\cD']$
generate the same subgroup of $\Br(K)$.

\qed

We have the following simple corollaries.

\begin{corollary}  \label{cor: fin-gen}
Let $K$ be a field as in Theorem \ref{th:main} and $L/K$ a finitely generated field extension.
Let also $\cD$ be a central division $L$-algebra. Then ${\bf gen}(\cD)$ consists of Brauer classes $[\cD']$ such that $[\cD]$ and $[\cD']$ generate the same subgroup of $\Br(L)$.

In particular, the genus of any division $L$-algebra of exponent 2 is trivial.
\end{corollary}

\noindent {\it Proof}. The field $L$ satisfies the same properties as the field $K$. That is, $L$
is a finite extension of a field which is a purely transcendental extension of infinite transcendence degree of some subfield.
\qed

\begin{corollary}  \label{cor: fin-gen}
Let $K$ be a field as in Theorem \ref{th:main}.
Then the genus of any division algebra of exponent 2 over the field of rational functions
$K(x)$ is trivial.
\end{corollary}


\section{Genus of simple algebraic groups of type $\mathrm{G}_2$} \label{sec:proof}


\noindent {\it Proof of Theorem \ref{th:G2}}.
Let $G$ and $H$ be simple algebraic groups of type $\mathrm{G}_2$ over $K$.
There exists a finitely generated field extension $E/F$ and simple algebraic groups $G'$ and $H'$ of type $\mathrm{G}_2$ over $E$ such that
$E \subset K$ and
$$
G \simeq G'\times_E K, H \simeq H' \times_E K.
$$
By Lemma \ref{l:finite subset},
there exists a field extension $M/E$ such that
$K/M$ is a purely transcendental field extension with infinite transcendence degree.

Let $S$ be a transcendence base of $K$ over $M$. Let also $s_1,\dots,s_6$ be some elements from $S$
and
$$
L:= M (S \setminus \{ s_1,\dots,s_6 \}).
$$
Then $K=L(s_1,\dots,s_6)$.

Moreover, $G \simeq (G'\times_E L)  \times_L K$ and $H \simeq (H'\times_E L) \times_L K$.
Assume that $G$ and $H$ are in the same genus.
Now using the same arguments as in the last paragraph of the proof of  Theorem 1.6 from \cite{ChRaRa24}, we obtain that
$G'\times_E L$ and $H'\times_E L$ are isomorphic over $L$. Hence $G$ and $H$ are also isomorphic.
\qed

\bigskip

{\bf Acknowledgements.} The author is grateful to A.S. Rapinchuk for the valuable comments and
suggestion to consider the genus problem for algebraic groups of type $\mathrm{G}_2$.



\end{document}